\documentclass{article}
\usepackage{amsfonts}
\usepackage{amssymb, amscd}

\newcommand{\RR}{{\cal R}}
\newcommand{\tr}{\mbox{tr}}

\newcommand{\uu}{{\mathfrak u}}

\newcommand{\R}{\mathbb{R}}
\newcommand{\Z}{\mathbb{Z}}

\newcommand{\C}{\mathbb{C}}

\newcommand{\pf}{\noindent{\bf Proof.} }

\newcommand{\qed}{{$\Box$}}
\newcommand{\W}{{\cal W}}

\newcommand{\mn}{\medskip\noindent}

\newcommand{\bn}{\bigskip\noindent}
\newcommand{\note}{\bn{\bf Note.} }

\newtheorem{cor}{Corollary}
\newtheorem{teo}{Theorem}

\newtheorem{lema}{Lemma}

\title{Orthogonal almost-complex structures of \\ minimal energy}
\author{Gil Bor, Luis Hern\'andez-Lamoneda\footnote{G.B and L.H.L acknowledge support from CONACyT grant 46274-F}, Marcos Salvai\footnote{M.S. acknowledges support from the following sources: \textsc{fonc}y\textsc{t}, Antorchas, \textsc{ciem\thinspace (conicet)} and
\textsc{sec}y\textsc{t\thinspace (unc)}.}}

\begin{document}
\maketitle

\begin{abstract}
In this article we apply a Bochner type formula  to show that on a
compact conformally flat riemannian manifold (or half-conformally
flat in dimension 4) certain types of orthogonal almost-complex
structures, if they exist, give the absolute minimum for the energy
functional. We give a few examples when such minimizers exist, and
in particular, we prove that the standard almost-complex structure
on the round $S^6$ gives the absolute minimum for the energy. We
also discuss the uniqueness of this  minimum and the extension of
these results to other orthogonal $G$-structures.
\end{abstract}
\section{Introduction}
Let $(M^{2n},g)$ be a Riemannian manifold. An orthogonal
almost-complex structure on $M$ is an automorphism of the tangent
bundle $J:TM\to TM$ which is orthogonal with respect to $g$ and
satisfies $J^2=-id_{TM}$. The combination  $(g,J)$ is also called an
``almost-hermitian structure'' on $M$.

Associated with such a structure is the K\"ahler form $\omega
=g(J\cdot ,\cdot)$ (or ``$J$ with its indices lowered by $g$'') and
the energy  $E(\omega)$, defined for a compact  manifold  by
$$
E(\omega)=\int_M{\| \nabla \omega\|^2}vol, \;
$$ where $\nabla\omega$ is the covariant derivative of $\omega$ with respect to the
Levi-Civita connection associated with $g$, and $vol$ is the volume
element  associated with $g$.\footnote{Thinking of $J$ as a section
of the twistor fibration over $M$ endowed with its natural
riemannian metric, this definition is equivalent to $E'(J)=\int_M{\|
dJ\|^2}$; i.e. $E'=aE+b$ where $a,b$ are a pair of constants
depending only on the dimension of $M$.}

A natural problem to consider in this context is that of  the
critical points of the energy functional, for a fixed $(M,g)$. In
particular,  one seeks orthogonal almost-complex structures $J$ of
minimal energy.

For $n=1$, i.e. $(M,g)$ an  oriented riemannian surface, $J$ is
unique (up to a sign) and $E(\omega)= 0$. For $n>1$, $E(\omega)\geq
0$  with equality if and only if  $\nabla\omega\equiv0$, which  is
the K\"ahler condition, i.e. $J$ is integrable and $\omega$ is
closed.

If $(M,g)$ does not admit a K\"ahler metric then we do not know in
general if a minimum occurs, let alone its value. However, in case
$(M,g)$ is {\em conformally flat} (or ``half-conformally-flat" for
$\dim M=4$), we are able to derive a useful sufficient condition for
the existence of an energy  minimizing $J$  and a formula for its
energy in terms of the total scalar curvature of $g$. This is the
main result of this paper (Theorem 1).


\bn The key ingredient for the proof is a general Bochner-type
formula for orthogonal $G$-structures previously published in [H]
and [BH]. Briefly, we consider the Gray-Hervella decomposition of
$\nabla \omega$, i.e. its decomposition  into the direct sum of four
$U_n$-irreducible components
$$\nabla\omega=\sum_{i=1}^4 (\nabla\omega)_i,$$
and the corresponding
$$E(\omega)=\sum_{i=1}^4 E_i(\omega),$$ where
$$E_i(\omega)=\int_M{\| (\nabla\omega)_i\|^2}vol,\quad i=1,\ldots, 4.$$ The Bochner-type formula of [BH] implies,
under the stated condition of conformal flatness on $(M,g)$, that
\begin{equation}2E_1(\omega)-E_2(\omega)+(n-1)E_4(\omega)=const., \end{equation}
where  ``$const.$'' is some  positive multiple of the total scalar
curvature of  $(M,g)$.

It follows immediately  from Formula (1)  that the vanishing of
certain components of $\nabla\omega$ implies that $J$ is an energy
minimizer. In particular, it follows  from formula (1) that {\em the
standard orthogonal almost-complex-structure on the 6-sphere $S^6$,
equipped with its standard (round) metric, is an energy  minimizer}
(see Theorem 2). Incidently, this result contradicts that of [W],
where  it is claimed  that this structure is not even a local
minimizer.

\bigskip

We notice that our formula (1) and its implications  is quite
similar to other variational problems of geometric origin, such as
the Yang-Mills equations in dimension 4 and harmonic maps between
K\"ahler manifolds.  In these problems, as in ours, there is a
natural decomposition of the ``energy'' into several components, and
one can show that a certain linear combination of these components
is identically constant. It follows that structures for which
certain components of the energy vanish are absolute minima. In this
way one sees that self-dual and anti-self-dual connections on
4-manifolds (the instantons) form the minima of the Yang-Mills
functional, and holomorphic or anti-holomorphic maps between
K\"ahler manifolds are the minima of the harmonic map energy (in
their homotopy class). Furthermore, in these cases, as in ours, the
condition for being a minimum is a first order differential equation
on the structure in question, whereas the Euler Lagrange equations
for general critical points of the energy are second-order PDE.

\bigskip

In the rest of this article, we first explain how to arrive at
equation (1) and the precise conditions under which it applies, and
then use formula (1) to give several examples of orthogonal
almost-complex structures that realize the absolute minimum of the
energy.
\bigskip

In the last section of the article we explain how to extend our
results to similar problems of ``$G$-structures with minimal
energy'', such as $G_2$ and $Spin_7$ structures.

\bn

{\bf Acknowledgments.} We thank V. Apostolov, M. Pontecorvo and S.
Simanca for helpful insights as well as some important references.

\section{A Bochner formula}
Let $(M^{2n},g,J)$ be an almost-hermitian manifold and $\omega
=g(J\cdot ,\cdot )$ its associated K\"ahler form. We give here a
brief review of the results of [BH] concerning such structure.

We use the abbreviated notation $\Lambda^k$ for the bundle of real
$k$ forms on $M$, $\Lambda^{p,q}$ for the bundle of complex forms of
type $(p,q)$ and $[[\Lambda^{p,q}]]$ for the real forms in
$\Lambda^{p,q}$.

In this notation, $\omega$ is a section of $[[\Lambda^{1,1}]]$, and
$\nabla\omega$ is a section of the bundle $\W:= \Lambda^1\otimes
[[\Lambda^{2,0}]]$. This bundle decomposes into four subbundles
$$\W=\W_1\oplus \ldots\oplus\W_4,$$
the so-called Gray-Hervella decomposition [GH], corresponding to the
decomposition into irreducibles of the $U_n$-representation which
gives rise to $\W$. We have

\begin{itemize}
\item $\W_1=[[\Lambda^{3,0}]]$;
\item $\W_2=$the real part of the image of
  $\left(
\Lambda^{1,0}\right)^{\otimes 3}$ under the Young symmetrizer
$(1-(23))(1 + (12))$;
\item  $\W_3=$ real part of the ``primitive'' part of $\Lambda^{1,2}$
  (kernel of the contraction in the first and second
  entries);
\item  $\W_4\cong\Lambda^1$,  given by the image of $\Lambda^1$, inside of
$\W$, of the adjoint of the contraction ${\cal
W}\subset\Lambda^1\otimes\Lambda^2\to\Lambda^1$.
\end{itemize}

\note  For $n=1$, $\W=0$; for $n=2$, $\W_1=\W_3=0$.

\mn

Corresponding to the decomposition of $\W$ is the decomposition of
$\nabla\omega$,
$$
\nabla\omega =(\nabla\omega)_1 +\ldots+(\nabla\omega)_4,
$$ i.e.  $(\nabla\omega)_i$ is a section of $\W_i$, $i=1,\ldots, 4$. It is important to notice that
the irreducible $U_n$-modules  giving rise to $\W_i$ are pairwise
non-isomorphic, hence  the $(\nabla\omega)_i$ are pairwise
orthogonal.

\medskip

The components $(\nabla\omega)_i$
  carry important geometric information about the
almost-complex structure; for example, $J$ is integrable iff
$(\nabla\omega)_1=(\nabla\omega)_2=0, $ i.e. $J\in \W_3\oplus \W_4$
(or  ``$J$ is of type $\W_3\oplus \W_4$''). The structure is
symplectic ($d\omega =0$), or ``almost Kahler'',  iff
$\nabla\omega\in\W_2$ and is  nearly-K\"ahler if
$\nabla\omega\in\W_1$, i.e. $\nabla\omega=d\omega$.

\medskip

Now when $M$ is compact, then corresponding to the decomposition of
$\nabla\omega$ is the decomposition of  the energy,
$$
E(\omega)=E_1(\omega)+E_2(\omega)+E_3(\omega)+E_4(\omega),
$$where
$$
E_i(\omega):=\int_M{\|(\nabla\omega)_i\|^2}.
$$

In [H] and [BH], via two different arguments, the following formula
for an arbitrary $(M^{2n},g,J)$ was obtained:
$$
2E_1(\omega)-E_2(\omega)+(n-1)E_4(\omega) = {1\over 2}\int_M{\tr
(\RR ,\uu_n^\perp )},$$ where  $\tr (\RR ,\uu_n^\perp )$ means the
trace of the $(\uu_n^\perp, \uu_n^\perp )$-block of the curvature
operator $\RR :\Lambda^2=\uu_n\oplus\uu_n^\perp\to
\uu_n\oplus\uu_n^\perp$.

\note  For $n=2$, since  only the $E_2$ and $E_4$ components exist,
i.e. $E(\omega)=E_2(\omega)+ E_4(\omega)$, the formula reduces to
$$
-E_2(\omega)+E_4(\omega) = {1\over 2}\int_M{\tr (\RR ,\uu_n^\perp
)}.$$

\mn

In general, the decomposition $\Lambda^2=\uu_n\oplus\uu_n^\perp$
depends on $J$, hence the same dependence occurs for the right-hand
side of the above formula. Nevertheless, when $(M,g)$ is conformally
flat (or half conformally flat in dimension 4), we have the
following:
\begin{lema} Let $(M^{2n},g)$ be a riemannian manifold with Weyl tensor $W$  and an orthogonal almost-complex structure $J$. If
\begin{itemize}
\item $n\geq 3$ and $W=0$  (i.e. $(M,g)$ is conformally flat), or
\item $n=2$ and  $W^+=0$ (i.e. $(M,g)$ is half-conformally flat, or anti-self-dual, using the orientation induced by $J$),
\end{itemize}
then
$$
\tr (\RR ,\uu_n^\perp )={2n-2\over 2n-1}s,
$$where $s$ is the scalar curvature.
\end{lema}

\mn {\bf Proof.} This is well known  (see for example, [G] or
[dRS]),  so we give here only  a sketch: by definition,  $\tr (\RR
,\uu_n^\perp )$ is a $U_n$-invariant functional on the space of
curvature type tensors. Representation theory tells us that there
are two linearly independent  such invariants, but restricted to the
space of curvature type tensors with vanishing  Weyl tensor $W$ (or
vanishing  $W^+$ in dimension 4) there is a unique $U_n$ invariant
(up to a constant). The exact value of the  constant may be
evaluated by computing it on any example (we used the real
hyperbolic $2n$-space). \qed
\begin{cor}Let $(M^{2n},g)$ be a compact Riemannian manifold such
that
\begin{itemize}
\item $n\geq 3$ and $(M,g)$ is conformally flat, or
\item $n=2$ and $(M,g)$ is anti-self-dual.
\end{itemize}
Then every almost-complex structure $J$, orthogonal with respect to
$g$, satisfies
$$
2E_1(\omega)-E_2(\omega)+(n-1)E_4(\omega) = C_g,
$$where $C_g$ is a constant depending only on the metric $g$;
in fact, $C_g={n-1\over 2n-1}\int_M{s}$, where $s$ is the scalar
curvature of $g$.
\end{cor}

\note For $n=2$, since $E_1=0$, the above formula reduces to
$$
-E_2(\omega)+E_4(\omega) = C_g.
$$

\mn

We now state our main result:

\begin{teo}Let $(M^{2n},g)$ be a compact Riemannian manifold such
that
\begin{itemize}
\item $n\geq 3$ and $(M,g)$ is conformally flat, or
\item $n=2$ and $(M,g)$ is anti-self-dual.
\end{itemize}
Then an orthogonal almost-complex structure $J_0$ on $M$ is an
energy  minimizer in each of the following 3 cases:
\begin{enumerate}
\item  $n=3$ and $J_0$ is of type $\W_1\oplus \W_4$.
\item $n=2$ or $n\geq 4$ and $J_0$ is of type $\W_4$.

\item $n$ is arbitrary and $J_0$ is of type ${\cal
W}_2$.
\end{enumerate}
Furthermore,
\begin{itemize}

\item $E(J_0)={1\over n-1}C_g$ in each of the first two cases,
$E(J_0)=-C_g$ for the third;

\item  if one of the above types of minimizers exists on $(M,g)$,
then any other minimizer is necessarily of the same type.
\end{itemize}
\end{teo}

\bn{\bf Remarks.}
\begin{enumerate}

\item For $n=2$, type $\W_4$  means $J_0$ is integrable.

\item For all $n\geq 2$, type ${\cal
W}_2$ means the associated  Kahler form $\omega_0$ is closed, i.e.
symplectic.

\item  Note that $\W_1$ is also of type
$\W_1\oplus \W_4$, hence for $n=3$, a $J_0$ of type $\W_1$ is a
minimizer. But the theorem does not exclude in this case the
existence of another minimizer of type $\W_1\oplus \W_4$ which is
not of type $\W_1$ (see the example of $S^6$ below).

\item  A conformal change of the metric only affects the $\W_4$ component of $\nabla\omega$ (see for example [GH]). Hence any minimizing $J_0$ of the first two types in the theorem  is an energy  minimizer with respect to all metrics in the conformal class of $g$. See more about this in section 3.5 below.

\item  There do exist riemannian  manifolds that do not
admit almost-complex structures of the indicated types, thus for
them it is not clear if minimizers exist or not (see the example
below of hyperbolic 4-manifolds).

\end{enumerate}

\bn{\bf Proof.}  Let  $J$ be an arbitrary  almost-complex structure
orthogonal with respect to $g$ and let $\omega,\omega_0$ be the
K\"ahler forms of $J,J_0$ (resp.).

\begin{enumerate}
\item If $n=3$ and $J_0$ is of type $\W_1\oplus {\cal
W}_4$, then
$$\begin{array}{rl}
  E(\omega)=&E_1(\omega)+E_2(\omega)+E_3(\omega)+E_4(\omega) \\ &\\
  \geq& E_1(\omega)-{1\over 2}E_2(\omega)+E_4(\omega)\\&\\
  =& C_g/2=E_1(\omega_0)-{1\over
2}E_2(\omega_0)+E_4(\omega_0) \\&\\
  =& E_1(\omega_0)+E_4(\omega_0)=E(\omega_0).
\end{array}
$$
\item If  $n=2$ or $n\geq 4$ and $J_0$ is of type ${\cal
W}_4$, then
$$
\begin{array}{rl}
  E(\omega)=&E_1(\omega)+E_2(\omega)+E_3(\omega)+E_4(\omega) \\ &\\
  \geq& {2\over n-1}E_1(\omega)-{1\over n-1}E_2(\omega)+E_4(\omega)\\&\\
  =& {C_g\over n-1}={2\over n-1}E_1(\omega_0)-{1\over
n-1}E_2(\omega_0)+E_4(\omega_0) \\&\\
  =& E_4(\omega_0)=E(\omega_0).
\end{array}
$$
\item Finally, if  $J_0$ is of type $\W_2$ (symplectic) then
$$
\begin{array}{rl}
  E(\omega)=&E_1(\omega)+E_2(\omega)+E_3(\omega)+E_4(\omega) \\ &\\
  \geq& -2E_1(\omega)+E_2(\omega)-(n-1)E_4(\omega)\\&\\
  =& -C_g=-2E_1(\omega_0)+E_2(\omega_0)-(n-1)E_4(\omega_0) \\&\\
  =& E_2(\omega_0)=E(\omega_0).
\end{array}
$$
\end{enumerate}

One can also read out easily from these calculations the exact value
of the minimal energy. For example, for $n=3$ and a $J_0$ of type
$\W_1\oplus {\cal W}_4$,
$$E(\omega_0)=E_1(\omega_0)+E_4(\omega_0)=
E_1(\omega_0)-{1\over 2}E_2(\omega_0)+E_4(\omega_0)=
 C_g/2.$$ The other cases are handled similarly.

 \mn

For the last statement of the theorem, notice (for example) that if
$n=3$ and a $J_0$ of type $\W_1\oplus {\cal W}_4$ exists, then if
$J$ is another minimizer we will have equality in the first
inequality above; thus, $E_2(\omega )=E_3(\omega)=0$ and $J$ is of
type $\W_1\oplus {\cal W}_4$. The other cases are handled similarly.
\mbox{\qed}

\section{Examples}

Here we  go through the 4 classes of minimal energy almost-hermitian
structures indicated in Theorem 1 and try to find examples in each
case.

\subsection{$S^6$}

The round 6-sphere has a natural compatible almost-complex
structure, $J_C$, given by Cayley cross-product in $\R^7$ (thought
of imaginary Cayley numbers). It is known that such structure is
nearly-K\"ahler (i.e. of type $\W_1$), and so, according to the
previous theorem, realizes the absolute minimum of the energy among
all almost-complex structures orthogonal with respect to the round
metric. In fact, due to the Remark 4 of the last section we can say
a little more:
\begin{teo} An almost complex structure on $S^6$ which is conformally equivalent
to the Cayley almost complex  structure is an energy  minimizer with
respect to any metric conformally equivalent to the round metric.
\end{teo}

\bn{\bf Remark.} In [CG] it has been shown that $J_C$ is  also a
{\em volume}  minimizer, among all  sections of the twistor
fibration.

\mn

Concerning the uniqueness of this minimizer we have the following
theorem of Friedrich [F]:
\begin{teo}  $J_C$ is the unique nearly-K\"ahler structure on the round
$S^6$, up to an isometry; i.e. for any almost complex structure $J$
on $S^6$ which is nearly-K\"ahler wrt the round metric there exists
an isometry $\phi\in SO_7$ such that $J=\phi^*J_C$.
\end{teo}

Hence $J_C$ is the unique (up to an isometry) energy minimizer on
the round $S^6$ of type $\W_1$. We know from theorem 1 that any
other minimizer should be   of type $\W_1\oplus\W_4$, but we do not
know if there is actually any one which is not conformally
equivalent  to $J_C$.


\subsection{Hopf manifolds}
The product metric on $S^{2n-1}\times S^1$ is conformally flat.
Moreover, this manifold admits many orthogonal complex structures
coming from actions of $\Z$ on $\C^n\setminus\{0\}$ by conformal
linear maps. Hence all of these structures are locally conformally
K\"ahler, thus of type $\W_4$ (see [V]); therefore they realize the
minimum for the energy.


\subsection{Symplectic manifolds}
It is well known that a conformally flat, K\"ahler manifold of
complex dimension at least three, must be flat [Be, 2.68]. Whereas
in dimension 2, one has also the product of two surfaces, one with
+1, the other -1, constant curvature. We don't know of any examples
of conformally flat, almost-K\"ahler (ie. symplectic), non-K\"ahler
compact manifolds.


\subsection{Anti-self-dual four-manifolds}
\subsubsection{Hermitian} The main result for
anti-self-dual-manifolds with orthogonal complex structure (the
anti-self-dual-hermitian, or ASDH, manifolds) is the following
theorem of Charles Boyer [B]:

\mn

{\it  if a compact ASDH-manifold has  even  first Betti number then
there is a conformal change of the metric that transform it into a
K\"ahler metric of zero scalar curvature.}


\bn

With odd first Betti number we have the following examples of LeBrun
[L]:

\mn

{\it There exist ASDH-metrics on the $k$-fold blow-ups
$$
(S^1\times S^3)\# \overline{\C P^2}\#\cdots \#\overline{\C P^2} .
$$}


\bn
\subsubsection{ASD Symplectic manifolds}
We recall Armstrong's deformation argument [A] to produce examples
of ASD  non-Kahler symplectic structures (we thank  V. Apostolov for
explaining this to us).


\bn

Start with a scalar flat Kahler (SFK) metric $(g_0, J_0)$ on a four
manifold. Such a metric can be shown to exist on certain complex
surfaces; for example,  on a blow-up of a generic ruled complex
surface, see [LKP]). Such a metric is ASD, since for a Kahler metric
ASD is equivalent to scalar flat.


\bn

 The idea is to show that on such a manifold there are  ASD deformations of the conformal class  $[g_0]$ admitting  symplectic structures which are not Kahler. To this end we consider two types
of deformations:
\begin{enumerate}
\item[(a)] SFK deformations of $(g_0,J_0)$.
\item[(b)]  ASD deformations of the conformal class $[g_0]$.
\end{enumerate}
These moduli spaces have been studied  by [LKP] and [KK]. Just by
comparing their dimensions we see that (b) is larger than (a) on our
manifold. So there are ASD classes $[g]$ arbitrarily close to
$[g_0]$ which are not Kahler. Let us see that such a $[g]$ admits a
compatible symplectic structure. Note first, that on a
four-manifold, the condition on a symplectic form $\omega$ to be
"admitted" by a conformal class $[g]$ is simply that  $\omega$ be
self-dual wrt $[g]$ (because this implies that $\omega$ is the
Kahler form associated with some orthogonal almost complex structure
and some metric in the conformal class of $g$). Now let $\eta$ be
the harmonic representative, wrt $g$, of the deRham cohomology class
of $\omega_0$, and let $\omega$ be the SD part of $\eta$. Then
$\omega$ is still harmonic, thus closed, and SD wrt $g$. If $g$ is
near $g_0$ then $\omega$ is near $\omega_0$ and is non-degenerate,
hence symplectic.

\subsection{Conformal change of the metric}

Consider an almost-hermitian manifold $(M,g,J)$, and a conformal
change of the metric: thus, set $g'=\lambda g$ with
$\lambda:M\to\R^+$. Note that $J$ is  orthogonal wrt $g'$ as well.

Denoting by $\omega '$ and $\nabla '$ the K\"ahler form and
Levi-Civita connection of the almost-hermitian manifold $(M,g',J)$,
we have the following relation (see [FFS] and [GH])
$$
\nabla '\omega '=\lambda\nabla\omega + \epsilon (\lambda)
$$where $\epsilon (\lambda)$ is a certain tensor in $\W_4$ depending on $\lambda$
and $d\lambda$ and such that $\epsilon (\lambda)\equiv 0$ if and
only if $\lambda$ is constant.

We see that there is no immediate relation that can be deduced
between the energies of $J$ with respect to $g$ and $g'$; in
particular, there's no obvious reason why a $J$ which has minimal
energy for $g$ should still be a minimizer for $g'$ -- even if
$(g,J)$ is a K\"ahler structure on $M$.

\bn


This motivates  the following question: {\it is there a compact
manifold $M$ and a conformal class of metrics $[g]$, such that $M$
is K\"ahler for two different (i.e. non-homothetic) metrics in
$[g]$?}

\bn

We do not know the answer to this question in general. However, from
Theorem 1 we know that if $(M,g,J)$ is a conformally flat
$\W_4$-manifold  (or ASDH-manifold in dimension 4, or ${\cal
W}_1\oplus \W_4$ conformally flat manifold, in dimension 6), then
$(M,g^\prime ,J)$ is still conformally flat (or ASD) and $J$ is of
the same type as for $g$. Thus $J$ will still be of minimum energy
for the new metric $g^\prime$. This leads to the following well
known result.

\begin{cor}Let $[g]$ be the conformal class of a metric on a compact
manifold. Assume that $[g]$ is ASD in dimension 4 or conformally
flat in higher dimensions. Then inside of $[g]$ there is at most one
K\"ahler metric (up to constant multiples).
\end{cor}


\bn

\pf Suppose $(g,J)$ is K\"ahler. Then $J$ has minimial energy
(namely 0) wrt $g$.  Let $g'=\lambda g$, for some non-constant
$\lambda:M\to \R^+$. Then $(g',J)$ is of type $\W_4$ and hence, by
Theorem 1, $J$ is  an energy minimizer wrt $g'$. Since $\lambda$ is
not constant then $\epsilon(\lambda)$ is not zero and so $(g',J)$
has positive energy. It follows that $g'$ cannot be Kahler, since
this would imply the existence of a $J'$ with zero energy wrt to
$g'$. \qed

\subsection{Hyperbolic 4-manifolds}
Let $(M^4,g)$ be a compact real hyperbolic 4-manifold. It is known
that many of these admit orthogonal almost complex structures (iff
its Euler characteristic is divisible by four), but cannot admit
neither a compatible complex structure ([G],[BH]) nor a compatible
symplectic structure ([OS]). We do not know whether these manifolds
have a minimizer or not.
\section{$G_2$ and $Spin_7$ structures}
A $G_2$-structure on a compact 7-manifold $M$ is given by a certain
3-form $\phi$. All $G_2$-structures giving the same metric $g$ are
given by sections of a fibre bundle with fiber $SO_7/G_2$ of
dimension $7$ (in fact, isometric to $\R P^7$). We can, again, talk
about the energy of a $G_2$-structure as the $L^2$-norm of
$\nabla\phi$.

In this case the formula is (see [BH]),
$$
6E_1+5E_7-E_{14}-E_{27}={2\over 3}\int_M{s}.
$$Thus any nearly parallel $G_2$-structure on $M^7$ minimizes energy
among all $G_2$-structures sharing the same metric. For example, the
standard $G_2$-structure on $S^7$ is minimal among all those
inducing the round metric.

Similarly, if $M^7$ admits a $\W_{14}\oplus \W_{27}$
$G_2$-structure, that structure will have minimal energy.

\bigskip

For $Spin_7$-structures the formula reads
$$
6E_8-E_{48}={1\over 6}\int_Ms.
$$Hence the natural $Spin_7$-structure on $S^7\times
S^1$  has minimal energy since it is of type $W_8$ [C].

\bigskip

\section{References.}

\mn [A] J. Armstrong, {\it Almost K\"ahler Geometry}, Ph.D.Thesis,
Oxford (1998).

\mn [B] C. Boyer, {\it Conformal Duality and Compact Complex
Surfaces}, Math. Ann. 274 (1986), 517-526.

\mn [Be] A. Besse, {\it Einstein manifolds} (Springer, Berlin,
1987).

\mn [BH] G. Bor and L. Hern\'andez-Lamoneda, {\it Bochner formulae
for orthogonal $G$-structures on compact manifolds}, Diff. Geom.
Appl. 15 (2001), 265-286.

\mn [C] F.M. Cabrera, {\it On Riemannian manifolda with
$Spin_7$-structure}, Publ. Math. Debrecen 46 (1995), 271-283.

\mn [CG] E. Calabi and H. Gluck, {\it What are the best
almost-complex structures on the $6$-sphere?}, Differential
geometry: geometry in mathematical physics and related topics (Los
Angeles, CA, 1990), 99--106, Proc. Sympos. Pure Math., 54, Part 2,
Amer. Math. Soc., Providence, RI, 1993.

\mn [dRS] H. del Rio and S. Simanca, {\it The Yamabe problem for
almost Hermitian manifolds}, J. Geom. Anal. 13 (2003), no. 1,
185--203.


\mn [F] T. Friedrich, {\it Nearly K\"ahler and nearly parallel
$G_2$-structures on spheres}, math.DG/0509146.

\mn [FFS] M. Falcitelli, A. Farinola and S. Salamon, {\it
Almost-Hermitian geometry}, Diff. Geom. Appl. 4 (1994), 259-282.

\mn [G] P. Gauduchon, {\it Complex structures on compact conformal
manifolds of negative type}, in: Complex Analysis and Geometry,
Lect. Notes Pure Appl. Math. 173 (1995), 201-212.

\mn [GH] A. Gray and L. Hervella, {\it The sixteen classes of almost
Hermitian manifolds and their linear invariants}, Ann. Mat. Pura
Appl. (4) 123 (1980), 35--58.

\mn [H] L. Hern'andez-Lamoneda, {\it Curvature vs. almost-Hermitian
structures}, Geom. Dedicata 79 (2000), 205-218.

\mn [K] M. Kapovich, {\it Conformally flat metrics on 4-manifolds},
J. Differential Geom. 66 (2004), no. 2, 289--301.

\mn [KK] A. King and D. Kotschick, {\it The deformation theory of
anti-self-dual conformal structures}, Math. Ann. 294 (1992), no. 4,
591--609.

\mn [L] C. LeBrun, {\it Anti-self-dual Hermitian metrics on blow-up
Hopf surfaces},  Math. Ann. 289 (1991), no. 3, 383--392.

\mn [LKP] C. LeBrun, J. Kim and M. Pontecorvo, {\it Scalar-flat
K\"ahler surfaces of all genera}, J. Reine Angew. Math. 486 (1997),
69--95.

\mn [OS] T. Oguro and K. Sekigawa, {\it Non-existence of almost
K\"ahler structure on hyperbolic spaces of dimension $2n(\geq 4)$},
Math. Ann. 300 (1994), no. 2, 317--329.

\mn [V] I. Vaisman, {\it On locally conformal almost K\"ahler
manifolds}, Israel J. Math. 24 (1976), 338-351.

\mn [W] C.M. Wood, {\it Instability of the nearly-K\"ahler
six-sphere}, J. Reine Angew. Math. 439 (1993), 205-212.

\end{document}